\documentclass[12pt]{article}

\usepackage{amsmath}

\usepackage{amssymb}

\usepackage{exscale}

\usepackage{amsthm}

\usepackage{epsfig}

\usepackage{verbatim}

\DeclareMathOperator{\disc}{disc}

\DeclareMathOperator{\conv}{conv}

\theoremstyle{plain}

\newtheorem{theorem}{Theorem}

\newtheorem{lemma}[theorem]{Lemma}

\newtheorem{corollary}[theorem]{Corollary}

\newcommand{\setmid}{\,|\,}

\newcommand{\N}{{\mathbb{N}}}

\newcommand{\R}{{\mathbb{R}}}

\newcommand{\EE}{{\mathcal{E}}}

\newcommand{\ee}{\varepsilon}

\begin{document}

\title{Balanced Partitions of Vector Sequences}


\author{Imre B\'ar\'any\thanks{R\'enyi Institute, Budapest,
POBox 127, 1364 Hungary, and Department of Mathematics,
University College London, Gower Street, London WC1E 6BT, United Kingdom.} \and Benjamin
  Doerr\thanks{Mathematisches Seminar, Bereich II,
    Christian-Albrechts-Universit\"at zu Kiel, 24098 Kiel, Germany.}}

\maketitle

\begin{abstract}
  Let $d, r \in \N$, $\|\cdot\|$ any norm on $\R^d$ and $B$ denote the unit ball
  with respect to this norm. We show that any sequence $v_1,v_2,\dots$ of
  vectors in $B$ can be partitioned into $r$ subsequences $V_1, \ldots, V_r$ in
  a balanced manner with respect to the partial sums: For all $n \in \N$, $\ell
  \le r$, we have $\|\sum_{i \le k, v_i \in V_\ell} v_i - \tfrac 1r \sum_{i \le
    k} v_i\| \le 2.0005 d$. A similar bound holds for partitioning sequences of
  vector sets. Both results extend an earlier one of B\'ar\'any and Grinberg
  (1981) to partitions in arbitrarily many classes.
\end{abstract}

\maketitle

\section{Introduction}

Let $d, N \in \N$.  Let $\|\cdot\|$ be any norm on $\R^d$ and $B = \{v \in \R^d
\setmid \|v\| \le 1\}$ its unit ball. In this paper, we give extensions of the
B\'ar\'any--Grinberg theorem to partitions into more than two classes. In its
most general version, this theorem states the following \cite{bg}.

\begin{theorem} \label{tselect}
  Let $V_1, \ldots, V_N \subseteq B$ such that $0 \in \conv(V_i)$ for all $i \in
  [N]$. Then there are $v_i \in V_i$ such that for all $n \in [N]$, \[\Big\|\sum_{i
  \in [n]} v_i\Big\| \le 2d.\]
\end{theorem}

The most interesting special case of Theorem~\ref{tselect} is that all $V_i$ are
of the form $V_i = \{v_i, -v_i\}$, cf. \cite{bs} as well. In this case,
Theorem~\ref{tselect} yields that for any sequence $v_1, \ldots, v_N$ of vectors
in $B$ there are signs $\ee_i \in \{-1,1\}$ such that $\|\sum_{i \in [n]} \ee_i
v_i\| \le 2d$ for all $n \in [N]$. In other words, there is a partition $[N] =
I_1 \dot\cup I_2$ such that $\|\sum_{i \in I_j \cap [n]} v_i - \tfrac 12 \sum_{i
  \in [n]} v_i\| \le d$ for all $n \in [N]$ and $j \in [2]$. This partitioning
version of the B\'ar\'any--Grinberg theorem was extended to partitions into $r >
2$ classes with error bound $(r-1) d$ in~\cite{wirmcol}. In the following
section, we show that the factor $(r-1)$ can be replaced by a constant.

In the third section of this paper, we show that if the stronger condition
$\sum_{v \in V_i} v = 0$ (instead of $0 \in \conv(V_i)$)
holds for all $i \in [N]$, then for each $i \in [N]$
there are $r$ distinct vectors $v_{i\ell} \in V_i$, $\ell \in [r]$, such that
$\|\sum_{i \in [n]} v_{i\ell}\| \le 5d$ holds for all $n \in [N]$ and $\ell \in
[r]$, where $r \le \max\{|V_i| \setmid i \in [N]\}$.

It is worth mentioning here that the results hold for all norms in $\R^d$.
This is due to the fact that proofs
use linear dependences among some vectors, with the norm playing very little
role. But most likely, much better bounds are valid for particular norms. For
instance, it is conjectured that for the $r=2$ and Euclidean norm case the best
bound is of order $\sqrt d$. (This was proved by Spencer~\cite{sp} when $N=O(d)$,
but the interesting case is when $N$ is arbitrary.)

In the proofs of both results below we invoke the recursive method
of~\cite{wirmcol}, which states, roughly speaking, that if there is a good
bound for (weighted) $2$--partitions, then there is a slightly worse bound for $r$--partitions.
In ~\cite{wirmcol} this is only worked out in the context of hypergraph
coloring (discrepancy problem), but the proofs easily reveal that all results
hold as well for balancing vectors. Thus we have the following:

\begin{theorem}\label{trec}
  Let $r \ge 2$ be an integer. Let $V$ be a set of vectors and $\EE$ be a set of
  subsets of $V$. Assume that for all integers $1 \le r_1 < r_0 \le r$ and all
  $V_0 \subseteq V$ there is a $V_1 \subseteq V_0$ such that for all $E \in
  \EE$, \[\Big\|\sum_{v \in V_1 \cap E} v - \tfrac {r_1}{r_0} \sum_{v \in V_0
    \cap E} v \Big\| \le K.\] Then there is a partition $V = V_1 \dot\cup \ldots
  \dot\cup V_r$ such that for all $\ell \in [r]$ and $E \in \EE$ we have
  \[\Big\|\sum_{v \in V_\ell \cap E} v - \tfrac 1r \sum_{v \in V \cap E} v
  \Big\| \le C(r) K,\]
  where $C(r)$ is an absolute constant satisfying $C(r) \le 2.0005$ for all
  $r \in \N$.
\end{theorem}

\section{Vector Partitioning}

Assume $V$ is a finite or infinite sequence of vectors $v_1,v_2,\dots$. We
introduce the (non-standard) notation $\sum _k V=\sum _{i=1}^k v_i$. Further,
for a subsequence $X$ of $V$ we define $\sum _k X=\sum _{i\leq k,\; v_i \in X}
v_i$.

\begin{theorem}
  For every sequence $V \subset B$, and for every integer $r\geq 2$, there is a
  partition of $V$ into $r$ subsequences $X_1,\dots,X_r$ such that for all $k$
  and $j$
  \[
  \sum_kX_j \in \frac 1r \sum_kV + C(r)dB.
  \]
\end{theorem}

\begin{proof}
  Assume $r=r_1+r_2$ (with positive integers $r_1,r_2$). We are going to
  construct a partition of $V$ into subsequences $Y_1$ and $Y_2$ such that for
  each $k$ and for $j=1,2$,
  \[
  \sum_kY_j \in \frac {r_j}{r} \sum_kV + dB.
  \]
  This implies the theorem via Theorem~\ref{trec}.

  For the construction of $Y_1,Y_2$ we use a modified version of the method of
  ``floating variables'' as given in~\cite{bg}.  Define
  $V_k=\{v_1,v_2,\dots,v_{k+d}\}$, $k=0, 1,2,\dots$. We are going to construct
  mappings $\beta _k\;: V_k \to [-r_1,r_2]$ and subsets $W_k \subset V_k$ with
  the following properties (for all $k$):
  \begin{enumerate}
  \item $\sum _{V_k}\beta _k(v)v=0$,
  \item $\beta _k(v) \in \{-r_1,r_2\}$ whenever $v \in W_k$,
  \item $|W_k|=k$ and $W_k \subset W_{k+1}$.
  \end{enumerate}
  The construction is by induction on $k$. For $k=0$,  $W_0 = \emptyset$
  and $\beta_0 = 0$ clearly suffice.  Now assume that $\beta _k$ and $W_k$ have
  been constructed and satisfy (i) to (iii).  The $d+1$ vectors in $V_{k+1}
  \setminus W_k$ are linearly dependent, so there are $\alpha(v) \in \R$ not all
  zero such that
  \[
  \sum _{V_{k+1} \setminus W_k} \alpha (v)v=0.
  \]
  Putting $\beta_k(v_{k+d+1})=0$, we have
  \[
  \sum_{W_k}\beta_k(v)v+\sum _{V_{k+1} \setminus W_k} (\beta_k(v)+t\alpha (v))v=0
  \]
  for all $t \in \R$. For $t=0$ all coefficients lie in $[-r_1,r_2]$. Hence for
  a suitable $t=t^*$, all coefficients still belong to $[-r_1,r_2]$, and
  $\beta_k(v)+t\alpha (v) \in \{-r_1,r_2\}$ for some $v=v^* \in V_{k+1}\setminus
  W_k$.  Set now $W_{k+1}=W_k \cup \{v^*\}$ and $\beta_{k+1}(v)=\beta_k(v)$, if
  $v \in W_k$, and $\beta_{k+1}(v)=\beta_k(v)+t^*\alpha (v)$, if $v \in
  V_{k+1}\setminus W_k$.  Now $W_{k+1}$ and $\beta_{k+1}$ satisfy the
  requirements. Moreover, $\beta_{k+1}(v)=\beta_k(v)$ for all $v \in W_k$.

  We now define the subsequences $Y_1$ and $Y_2$. Put $v_i$ into $Y_1$ if
  $v_i \in W_k$ and $\beta_k(v_i)=r_2$ for some $k$, and put $v_i$ into $Y_2$ if
  $v_i \in W_k$ and $\beta_k(v_i)=-r_1$ for some $k$. As
  $\beta_k(v)=\beta_{k+1}(v)$ once $v \in W_k$, this definition is correct for
  all vectors that appear in some $W_k$.  The remaining (at most $d$) vectors
  can be put into $Y_1$ or $Y_2$ in any way.  Set $\gamma (v)=r_2$, if $v \in
  Y_1$, and $\gamma (v)=-r_1$, if $v \in Y_2$.

  Clearly, $r_2\sum_kY_1 -r_1\sum_kY_2 \in rdB$ for all $k\leq d$. For $k>d$ we
  have, with $k=h+d$,
  \begin{eqnarray*}
    r_2\sum_kY_1 -r_1\sum_kY_2
    &=& \sum_{V_h}\gamma (v)v= \sum_{V_h}\gamma (v)v - \sum_{V_h}\beta_h(v)v \\
    &=& \sum_{V_h}(\gamma (v)- \beta_h(v))v=
    \sum_{V_h \setminus W_h}(\gamma(v)-\beta_h(v))v.
  \end{eqnarray*}
  The last sum contains at most $d$ non-zero terms, each having norm at
  most~$r$. Thus
  \[
  r_2\sum_kY_1 -r_1\sum_kY_2 \in rdB
  \]
  for every $k$. Adding this to the trivial equation $r_1\sum_kY_1 +r_1\sum_kY_2
  =r_1\sum_kV$ (expressing that $Y_1,Y_2$ form a partition of $V$), we obtain
  \[
  \sum_kY_1 \in \frac {r_1}{r}\sum_kV + dB
  \]
  for every $k$.
\end{proof}

\section{Vector Selection}

Let now $V_1, \ldots, V_{N}$ be a sequence of finite subsets of $B$ such that
$|V_i| \ge r$ for all $i \in [N]$. An $r$--selection of $(V_i)$ is a mapping
$\chi: [N] \times [r] \to \R^d$ such that $\chi(i,[r])$ is a $r$--element subset
of $V_i$ for all $i \in [N]$. For such a $\chi$, we define its discrepancy with
respect to $(V_i)$ by
\[\disc(\chi,(V_i)_{i \in [N]}) = \max_{n \in[N]} \max_{\ell \in [r]} \Big\|\sum_{i \in
  [n]} \Big(\chi(i,\ell) - \frac 1 {|V_i|} \sum_{v \in V_i} v\Big)\Big\|.\]

\begin{theorem}\label{tselectnew}
  There is an $r$--selection with discrepancy at most $5 d$.
\end{theorem}

We mention that this theorem also holds for infinite sequences of finite subsets
of~$B$.

To prove the theorem, we apply the following lemma twice.

\begin{lemma}
  Let $r \in \N$, $r\ge 2$. Let $V_1, \ldots, V_N \subseteq B$ such that $|V_i|
  \ge r$ for all $i \in [N]$. Then for all $k \in [r]$ there are $U_i \subseteq
  V_i$ such that $|U_i| = k$ and $\max_{n \in [N]} \|\sum_{i \in [n]} (\sum_{v
    \in U_i} v - \tfrac k {|V_i|} \sum_{v \in V_i} v)\| \le 2d$.
\end{lemma}

\begin{proof}
  We give an algorithm for the construction the $U_i$. For each $i \in [N]$, $v
  \in V_i$ put $x_{iv} = \tfrac k {|V_i|}$. We iteratively change these numbers
  to zeros and ones in such a way that $U_i := \{v \in V_i \setmid x_{iv}=1\}$ gives
  the desired solution. For the start let $n = 1$. What we do is the following:
  View those $x_{iv}$ such that $x_{iv} \notin \{0,1\}$ and $i \le n$ as
  variables. If there is exactly one solution to the linear system
  \begin{eqnarray}
    &&\sum_{i \in [N]} \sum_{v \in V_i} (x_{iv} - \tfrac k {|V_i|}) v = 0 \label{eq1}\\
    &&\sum_{v \in [|V_i|]} x_{iv} = k, \quad i \in [N] \label{eq2}\\
    &&x_{iv} \in [0,1], i \in [N], v \in V_i,\nonumber
  \end{eqnarray}
  then increase $n$ by one and try again. Otherwise our existing solution may be
  changed in such a way that at least one more variable $x_{iv}$ becomes $0$ or
  $1$. If $n$ reaches $N$ and no solution can be found, then stop and change the
  remaining non-integral values of $x_{iv}$ to $0$ or $1$ in such a way that
  (\ref{eq2}) is still fulfilled.

  Assume that in some step of this iteration no solution can be found. Then
  there are at least as many constraints containing variables as there are
  variables. Let $q$ be the number of constraints of type (\ref{eq2}) that
  contain a variable. Then the total number of constraints containing variables
  is at most $d+q$, and the number of variables is at least $2q$. Hence $q \le
  d$ holds if no non-trivial solution can be found, and at most $q+d \le 2d$ of
  the $x_{iv}$, $i \le n$, are not in $\{0,1\}$. Denote the set of these pairs
  $(i,v)$ by $I$. Since the remaining $x_{iv}$, $i \le n$, are not changed
  anymore, our final solution $\tilde x$ satisfies
  \begin{eqnarray*}
    \sum_{i \in [n]} \sum_{v \in V_i} (\tilde x_{iv} - \tfrac k {|V_i|})
      v
          &=& \sum_{i \in [n]} \sum_{v \in V_i} (x_{iv} - \tfrac k {|V_i|})
    v + \sum_{(i,v) \in I} (\tilde x_{iv} - x_{iv}) v \\
    &=& \sum_{(i,v) \in I} (\tilde x_{iv} - x_{iv}) v.
   \end{eqnarray*}
   Since $|I| \le 2d$, we conclude $\|\sum_{i \in [n]} \sum_{v \in V_i} (\tilde
   x_{iv} - \tfrac k {|V_i|}) v\|\le 2d$ for all $n \in [N]$. Since $\tilde
   x_{iv} \in \{0,1\}$, putting $U_i := \{v \in V_i \setmid \tilde x_{iv}=1\}$
   gives the desired solution.
\end{proof}

\begin{proof}[Proof of the theorem]
  Let us assume first that $|V_i| = r$ for all $i \in [N]$. Then, by the above
  lemma, for all integers $r_1, r_2$ such that $r = r_1 + r_2$ there are
  $U_i^{(1)} \dot\cup U_i^{(2)} = V_i$ such that $|U_i^{(j)}| = r_j$ and
  $\|\sum_{i \in [n]}(\sum_{v \in U_i^{(j)}} v - \tfrac{r_j}{r} \sum_{v \in V_i}
  v) \| \le 2d$. Hence from Theorem~\ref{trec}, we obtain an $r$--selection
  (actually an $r$--partition) of $(V_i)$ such that
 \[\Big\|\sum_{i \in [n]} \Big(\chi(i,\ell) - \tfrac 1r \sum_{v \in V_i} v\Big)
  \Big\| \le 2 \, C(r) d\] for all $n \in [N]$, $\ell \in [r]$.

  If $|V_i| > r$ for some $i$, apply the lemma (with $k = r$) to obtain $\tilde
  V_i \subseteq V_i$ such that $|\tilde V_i| = r$ and $\|\sum_{i \in
   [n]}(\sum_{v \in \tilde V_i} v - \tfrac{r}{|V_i|} \sum_{v \in V_i} v) \| \le
  2d$. By the above, there is a $r$--selection for $(\tilde V_i)$ such that
  \[\Big\|\sum_{i \in [n]} \Big(\chi(i,\ell) - \tfrac 1r \sum_{v \in \tilde V_i} v\Big) \Big\|
  \le 2 \, C(r) d\] for all $n \in [N]$, $\ell \in [r]$. Note that, trivially,
  $\chi$ is also an $r$--selection for $(V_i)$. It satisfies
  \begin{eqnarray*}
    \lefteqn{\Big\|\sum_{i \in [n]} \Big(\chi(i,\ell) - \tfrac 1{|V_i|} \sum_{v \in
      V_i} v\Big)   \Big\|}\\
    &\le & \Big\|\sum_{i \in [n]} \Big(\chi(i,\ell) - \tfrac 1r \sum_{v \in  \tilde V_i} v\Big)
    \Big\| + \Big\|\sum_{i \in [n]} \Big(\tfrac 1r \sum_{v \in  \tilde V_i} v - \tfrac
    1{|V_i|} \sum_{v \in  V_i} v\Big)\Big\|\\
    &\le& 2 \, C(r) d + \tfrac 1r 2d
  \end{eqnarray*}
  for all $n \in [N]$ and $\ell \in [r]$. By noting that $C(2) = 1$ and $C(r)
  \le 2.0005$ for all $r \in \N$, we obtain the constant
  of~$5$.
\end{proof}

We may remark that a closer inspection of $C(r)$ for small $r$ yields better
constants. For example, easy calculations by hand or Lemma~3.5 in~\cite{wirmcol}
show that $C(r) + \tfrac 1r \le 2.1$ for $r \le 10$ (for $r=7$ observe that
$C(7) \le \max\{\tfrac 13 + C(3),\tfrac 14 + C(4)\}$). Hence the bound $C(r) \le
2.0005$ implies $C(r) + \tfrac 1r \le 2.1$ for all $r \in \N$, leading to a
constant of $4.2$ instead of $5$.

The following is an immediate consequence of Theorem~\ref{tselectnew}.

\begin{corollary}
  Let $r, N \in \N$. For $i \in [N]$ let $V_i \subseteq B$ such that $\sum_{v
    \in V_i} v = 0$ and $|V_i| \ge k$. Then there is a $k$--selection of $(V_i)$
  such that
  \[\Big\|\sum_{i \in [n]} \chi(i,\ell) \Big\| \le 5d\] for all $n \in [N]$ and $\ell
  \in [r]$.
\end{corollary}

This answers a question of Emo Welzl concerning multi-class extensions of
Theorem~\ref{tselect} posed at the Oberwolfach Seminar on ``Discrepancy Theory
and its Applications'' in March 2004. It is clear that the stronger assumption
$\sum_{v \in V_i} v = 0$ is necessary. Already for $d=1$ and $r = 2$, the
sequence $V_i = \{-\tfrac 12, 1\}$ shows that $0 \in \conv(V_i)$ does not
suffice.

\subsection*{Acknowledgment}

We thank the organizers of the Oberwolfach Seminar on ``Discrepancy Theory
and its Applications'' (March 2004) as well as the Oberwolfach crew for
providing us with surroundings that resulted in this paper.

The first named author is grateful to Microsoft Research, (in Redmond, WA) as
part of the research reported here was carried out on a very pleasant and
fruitful visit there. For the same nice reason, the second author would like
to thank Joel Spencer and the Courant Institute of Mathematical Sciences (New
York City).


\begin{thebibliography}{BG81}
\bibitem[BG81]{bg}
I.~B\'ar\'any and V.~S. Grinberg.
\newblock On some combinatorial questions in finite-dimensional spaces.
\newblock {\em Linear Algebra Appl.}, 41:1--9, 1981.



\bibitem[BS95]{bs}
J.~Beck and V.~T. S\'os.
\newblock Discrepancy theory.
\newblock In R.~Graham, M.~Gr\"otschel, and L.~Lov\'asz, editors, {\em Handbook
  of {C}ombinatorics}, pages 1405--1446. Elsevier, 1995.



\bibitem[DS03]{wirmcol}
B.~Doerr and A.~Srivastav.
\newblock Multicolour discrepancies.
\newblock {\em Combinatorics, Probability and Computing}, 12:365--399, 2003.

\bibitem[Sp86]{sp}
J.~Spencer.
\newblock Balancing vectors in the max norm.
\newblock {\em Combinatorica}, 6:55--65, 1986.

\end{thebibliography}
\end{document}